\magnification 1095
\input plainenc
\input amssym
\input xepsf
\fontencoding{T2A}
\inputencoding{utf-8}
\tolerance 4000
\relpenalty 10000
\binoppenalty 10000
\parindent 1.5em

\hsize 17truecm
\vsize 24truecm
\hoffset -0.25truecm
\voffset -0.5truecm

\font\TITLE labx1440
\font\tenrm larm1000
\font\cmtenrm cmr10
\font\tenit lati1000
\font\tenbf labx1000
\font\teni cmmi10 \skewchar\teni '177
\font\tensy cmsy10 \skewchar\tensy '60
\font\tenex cmex10
\font\teneufm eufm10
\font\eightrm larm0800
\font\cmeightrm cmr8
\font\eightit lati0800
\font\eightbf labx0800
\font\eighti cmmi8 \skewchar\eighti '177
\font\eightsy cmsy8 \skewchar\eightsy '60
\font\eightex cmex8
\font\eighteufm eufm8

\font\cmsixrm cmr6

\font\sixbf labx0600
\font\sixi cmmi6 \skewchar\sixi '177
\font\sixsy cmsy6 \skewchar\sixsy '60
\font\sixeufm eufm6

\font\cmfiverm cmr5

\font\fivebf labx0500
\font\fivei cmmi5 \skewchar\fivei '177
\font\fivesy cmsy5 \skewchar\fivesy '60
\font\fiveeufm eufm5
\font\tencmmib cmmib10 \skewchar\tencmmib '177
\font\eightcmmib cmmib8 \skewchar\eightcmmib '177
\font\sevencmmib cmmib7 \skewchar\sevencmmib '177
\font\sixcmmib cmmib6 \skewchar\sixcmmib '177
\font\fivecmmib cmmib5 \skewchar\fivecmmib '177
\newfam\cmmibfam
\textfont\cmmibfam\tencmmib \scriptfont\cmmibfam\sevencmmib
\scriptscriptfont\cmmibfam\fivecmmib
\def\tenpoint{\def\rm{\fam0\tenrm}\def\it{\fam\itfam\tenit}%
	\def\bf{\fam\bffam\tenbf}
	\textfont0\cmtenrm \scriptfont0\cmsevenrm \scriptscriptfont0\cmfiverm
  	\textfont1\teni \scriptfont1\seveni \scriptscriptfont1\fivei
  	\textfont2\tensy \scriptfont2\sevensy \scriptscriptfont2\fivesy
  	\textfont3\tenex \scriptfont3\tenex \scriptscriptfont3\tenex
  	\textfont\itfam\tenit
	\textfont\bffam\tenbf \scriptfont\bffam\sevenbf
	\scriptscriptfont\bffam\fivebf
	\textfont\eufmfam\teneufm \scriptfont\eufmfam\seveneufm
	\scriptscriptfont\eufmfam\fiveeufm
	\textfont\cmmibfam\tencmmib \scriptfont\cmmibfam\sevencmmib
	\scriptscriptfont\cmmibfam\fivecmmib
	\normalbaselineskip 12pt
	\setbox\strutbox\hbox{\vrule height8.5pt depth3.5pt width0pt}%
	\normalbaselines\rm}
\def\eightpoint{\def\rm{\fam 0\eightrm}\def\it{\fam\itfam\eightit}%
	\def\bf{\fam\bffam\eightbf}%
	\textfont0\cmeightrm \scriptfont0\cmsixrm \scriptscriptfont0\cmfiverm
	\textfont1\eighti \scriptfont1\sixi \scriptscriptfont1\fivei
	\textfont2\eightsy \scriptfont2\sixsy \scriptscriptfont2\fivesy
	\textfont3\eightex \scriptfont3\eightex \scriptscriptfont3\eightex
	\textfont\itfam\eightit
	\textfont\bffam\eightbf \scriptfont\bffam\sixbf
	\scriptscriptfont\bffam\fivebf
	\textfont\eufmfam\eighteufm \scriptfont\eufmfam\sixeufm
	\scriptscriptfont\eufmfam\fiveeufm
	\textfont\cmmibfam\eightcmmib \scriptfont\cmmibfam\sixcmmib
	\scriptscriptfont\cmmibfam\fivecmmib
	\normalbaselineskip 11pt
	\abovedisplayskip 5pt
	\belowdisplayskip 5pt
	\setbox\strutbox\hbox{\vrule height7pt depth2pt width0pt}%
	\normalbaselines\rm
}

\def\No{\char 157}
\def\empty{}

\catcode`\@ 11
\catcode`\" 13
\def"#1{\ifx#1<\char 190\relax\else\ifx#1>\char 191\relax\else #1\fi\fi}

\def\newl@bel#1#2{\expandafter\def\csname l@#1\endcsname{#2}}
\openin 11\jobname .aux
\ifeof 11
	\closein 11\relax
\else
	\closein 11
	\input \jobname .aux
	\relax
\fi

\newcount\c@section
\newcount\c@subsection
\newcount\c@subsubsection
\newcount\c@equation
\newcount\c@bibl
\newcount\c@enum
\c@section 0
\c@subsection 0
\c@subsubsection 0
\c@equation 0
\c@bibl 0
\newdimen\d@enum
\d@enum=0pt
\def\lab@l{}
\def\label#1{\immediate\write 11{\string\newl@bel{#1}{\lab@l}}%
	\ifhmode\unskip\fi}
\def\eqlabel#1{\rlap{$(\equation)$}\label{#1}}

\def\section#1{\global\advance\c@section 1
	{\par\vskip 3ex plus 0.5ex minus 0.1ex
	\rightskip 0pt plus 1fill\leftskip 0pt plus 1fill\noindent
	{\bf\S\thinspace\number\c@section .~#1}\par\penalty 25000%
	\vskip 1ex plus 0.25ex}
	\gdef\lab@l{\number\c@section.}
	\c@subsection 0
	\c@subsubsection 0
	\c@equation 0
}
\def\subsection{\global\advance\c@subsection 1
	\par\vskip 1ex plus 0.1ex minus 0.05ex{\bf\number\c@subsection. }%
	\gdef\lab@l{\number\c@section.\number\c@subsection}%
	\c@subsubsection 0\c@equation 0%
}
\def\subsubsection{\global\advance\c@subsubsection 1
	\par\vskip 1ex plus 0.1ex minus 0.05ex%
	{\bf\number\c@subsection.\number\c@subsubsection. }%
	\gdef\lab@l{\number\c@section.\number\c@subsection.%
		\number\c@subsubsection}%
}
\def\equation{\global\advance\c@equation 1
	\gdef\lab@l{\number\c@section.\number\c@subsection.%
	\number\c@equation}{\rm\number\c@equation}
}
\def\bibitem#1{\global\advance\c@bibl 1
	[\number\c@bibl]%
	\gdef\lab@l{\number\c@bibl}\label{#1}
}
\def\ref@ref#1.#2:{\def\REF@{#2}\ifx\REF@\empty{\S\thinspace#1}%
	\else\ifnum #1=\c@section {#2}\else {\S\thinspace#1.#2}\fi\fi
}
\def\ref@eqref#1.#2.#3:{\ifnum #1=\c@section\ifnum #2=\c@subsection
	{(#3)}\else{#2\thinspace(#3)}\fi\else{\S\thinspace#1.#2\thinspace(#3)}\fi
}
\def\ref#1{\expandafter\ifx\csname l@#1\endcsname\relax
	{\bf ??}\else\edef\mur@{\csname l@#1\endcsname :}%
	{\expandafter\ref@ref\mur@}\fi
}
\def\eqref#1{\expandafter\ifx\csname l@#1\endcsname\relax
	{(\bf ??)}\else\edef\mur@{\csname l@#1\endcsname :}%
	{\expandafter\ref@eqref\mur@}\fi
}
\def\cite#1{\expandafter\ifx\csname l@#1\endcsname\relax
	{\bf ??}\else\hbox{\bf\csname l@#1\endcsname}\fi
}

\def\im{\mathop{\rm im}}
\def\supp{\mathop{\rm supp}}
\def\Wo{{\mathpalette\Wo@{}}W}
\def\Wo@#1{\setbox0\hbox{$#1 W$}\dimen@\ht0\dimen@ii\wd0\raise0.65\dimen@%
\rlap{\kern0.35\dimen@ii$#1{}^\circ$}}

\catcode`\"=12
\def\bolddelta{\mathchar"0\hexnumber@\cmmibfam 0E}
\catcode`\"=13
\long\def\enumerate#1#2{%
	\setbox0\hbox{$#1^{\circ}.\ $}\d@enum=\wd0\global\advance\d@enum 2pt\c@enum=0%
	{\def\item{\global\advance\c@enum 1\par\hskip 0pt%
	\hbox to \d@enum{$\number\c@enum^{\circ}$.\hss}}%
	\par\smallskip #2\par\smallskip}
}
\def\proof{\par\medskip{\rm Д$\,$о$\,$к$\,$а$\,$з$\,$а$\,$т$\,$е$\,$л$\,$ь%
	$\,$с$\,$т$\,$в$\,$о.}\ }
\def\endproof{{\parfillskip 0pt\hfill$\square$\par}\medskip}

\catcode`\@ 12
\immediate\openout 11\jobname.aux


\frenchspacing\rm
\leftline{УДК~517.984}\vskip 5pt
{\TITLE\rightskip 0pt plus 1fill\leftskip 0pt plus 1fill\noindent
Об одном классе сингулярных задач Штурма--Лиувилля\par
\vskip 2ex\noindent\rm А.$\,$А.~Владимиров\footnote{}{\eightrm Работа поддержана
РФФИ, грант \No~13-01-00705.}\par}
\vskip 0.25cm
$$
	\vbox{\hsize 0.75\hsize\leftskip 0cm\rightskip 0cm
	\eightpoint\rm
	{\bf Аннотация:\/} Даётся определение самосопряжённых граничных задач,
	отвеча\-ющих формальному дифференциальному выражению
	$$
		-(y'/r)'+qy=\lambda py,
	$$
	где $r$~--- конечная борелевская мера, а $q$ и $p$~--- обобщённые функции
	из некоторого связанного с мерой $r$ класса. Описываются отвечающие такому
	пониманию неог\-раниченные самосопряжённые операторы в функциональных
	пространствах. В случае, когда $q=0$, а $r$ и $p$ самоподобны, указываются
	асимптотики спектра.
	}
$$

\vskip 0.5cm
\section{Общие конструкции}\label{par:2}
\subsection\label{pt:1.1}
Пусть $r\in W_2^{-1}[0,1]$~--- некоторая отличная от тождественного нуля
неотрица\-тель\-ная обобщённая функция, допускающая представление
$$
	\langle r,y\rangle\equiv-\int_0^1 R\overline{y'}\,dx+R(1)\overline{y(1)}-
		R(0)\overline{y(0)},
$$
где заданная почти всюду функция $R\colon[0,1]\to [R(0),R(1)]$ неубывает.
Пусть также $U$~--- некоторая унитарная комплексная матрица размера $2\times 2$.
Введём в рассмотрение гильбертовы пространства
$$
	\displaylines{W_{2,U}^1[R(0),R(1)]\rightleftharpoons
		\left\{u\in W_2^1[R(0),R(1)]\;:\;u^\wedge\in\im (U-1)\right\},\cr
	\frak D_{R,U}\rightleftharpoons\left\{\pmatrix{y\cr u}\in L_2[0,1]\times
		W_{2,U}^1[R(0),R(1)]\;:\; y=u\circ R\right\},\qquad
	\left\|\pmatrix{y\cr u}\right\|_{\frak D_{R,U}}\rightleftharpoons
		\|u\|_{W_2^1[R(0),R(1)]},}
$$
где через $u^\wedge$, следуя [\cite{RH:2001}, \hbox{(7.50)}], обозначен
связанный с функцией $u\in W_2^1[R(0),R(1)]$ вектор граничных значений
$$
	u^\wedge\rightleftharpoons\pmatrix{u(R(0))\cr u(R(1))}.
$$
Запись $y\in\frak D_{R,U}$ для функции $y\in L_2[0,1]$ мы будем далее понимать
как сокращённую запись условия
$$
	(\exists u\in W_{2,U}^1[R(0),R(1)])\qquad \pmatrix{y\cr u}\in\frak D_{R,U}.
$$
В случае непрерывности функции $R$ функция $u\in W_{2,U}^1[R(0),R(1)]$ однозначно
восстанавливается по заданной функции $y=u\circ R\in L_2[0,1]$. В общем случае
это уже не так. Тем не менее, там, где это не приводит к недоразумениям, мы будем
использовать запись $\|y\|_{\frak D_{R,U}}$ в качестве условного обозначения
величины $\|u\|_{W_2^1[R(0),R(1)]}$ и в общей ситуации.

Через $S_{R,U}\colon W_{2,U}^1[R(0),R(1)]\to\frak D_{R,U}$ мы будем обозначать
естественную изометриче\-скую биекцию пространств $W_{2,U}^1[R(0),R(1)]$
и $\frak D_{R,U}$.

\subsection
Может быть указана постоянная $C_R>0$, независимо от выбора функции
$u\in W_2^1[R(0),R(1)]$ удовлетворяющая оценке
$$
	\|u\|_{C[R(0),R(1)]}\leqslant C_R\cdot\|u\|_{W_2^1[R(0),R(1)]}.
	\leqno(\equation)
$$\label{eq:CD}%
Соответственно, неравенства
$$
	\eqalign{\left|\int_0^1 fy\cdot\overline{z}\,dx\right|&\leqslant
			\|f\|_{L_1[0,1]}\cdot\|y\|_{L_{\infty}[0,1]}\,
			\|z\|_{L_{\infty}[0,1]}\cr
		&\leqslant C_R^2\,\|f\|_{L_1[0,1]}\cdot
			\|y\|_{\frak D_{R,U}}\,\|z\|_{\frak D_{R,U}}}
$$
означают, что с каждой функцией $f\in L_1[0,1]$ может быть связан ограниченный
оператор умножения, отображающий пространство $\frak D_{R,U}$ в двойственное к нему
пространство $\frak D^*_{R,U}$. Полная непрерывность вложения пространства
$W_2^1[R(0),R(1)]$ в пространство $C[R(0),R(1)]$ гарантирует при этом полную
непрерывность указанного мультипликатора.

Через $\frak M_{R,U}$ мы далее будем обозначать пространство мультипликаторов,
получаемое замыканием относительно равномерной операторной топологии класса
отображающих $\frak D_{R,U}$ в $\frak D^*_{R,U}$ операторов умножения
на суммируемые функции. В свете вышесказанного, любой мультипликатор $f\in
\frak M_{R,U}$ является вполне непрерывным.

\subsubsection\label{prop:1.1}
{\it Всякая заданная на отрезке $[0,1]$ конечная борелевская мера $f$, не имеющая
с мерой $r$ общих атомов, определяет некоторый элемент пространства $\frak M_{R,U}$.
}%
\proof
Зафиксируем произвольное значение $\varepsilon>0$ и поставим ему в соот\-ветствие
такое разбиение
$$
	0=\zeta_0<\zeta_1<\zeta_2<\ldots<\zeta_n<\zeta_{n+1}=1
$$
интервала $(0,1)$ точками непрерывности меры $f+r$, что при любом выборе индекса
$k\in\overline{0,n}$ выполняется хотя бы одно из неравенств
$$
	\int_{\zeta_k}^{\zeta_{k+1}}f\,dx<\varepsilon^3\quad\hbox{или}\quad
	R\bigr|_{\zeta_k}^{\zeta_{k+1}}<\varepsilon^2.\leqno(\equation)
$$\label{eq:**}%
Разобьём систему всевозможных индексов $k\in\overline{0,n}$ на две подсистемы
$\Gamma$ и $\bar\Gamma$ таким образом, чтобы приращение функции
$R$ на любом отрезке $[\zeta_k,\zeta_{k+1}]$ при $k\in\Gamma$ минорировало
$4\varepsilon^2$, а при $k\in\bar\Gamma$~--- мажорировало $\varepsilon^2$.
Число элементов набора $\bar\Gamma$ при этом с очевидностью не может превосходить
величину $[R(1)-R(0)]\,\varepsilon^{-2}$. Наконец, введём в рассмотрение
ступенчатую функцию $f_{\varepsilon}$ вида
$$
	f_{\varepsilon}\rightleftharpoons\sum\limits_{k\in\Gamma}
		{\chi_{[\zeta_k,\zeta_{k+1}]}\over\zeta_{k+1}-\zeta_k}
		\int_{\zeta_k}^{\zeta_{k+1}}f\,dx,
$$
где через $\chi_\Lambda$ обозначен индикатор промежутка $\Lambda$. Проводимые
с учётом соотношений \eqref{eq:CD}, \eqref{eq:**} и интегральной теоремы о среднем
несложные вычисления показывают, что для любой точки $u$ единичного шара
пространства $W_2^1[R(0),R(1)]$ существуют промежуточные значения $\xi_k\in
[R(\zeta_k),R(\zeta_{k+1})]$, удовлетворяющие оценкам
$$
	\eqalign{\left|\int_0^1 (f-f_\varepsilon)\cdot|u\circ R|^2\,dx
		\right|&\leqslant\sum_{k\in\bar\Gamma}\int_{\zeta_k}^{\zeta_{k+1}}
		\kern -1mm f\cdot |u\circ R|^2\,dx+\sum\limits_{k\in\Gamma}
		\int_{\zeta_k}^{\zeta_{k+1}}\kern -1mm f\cdot\biggl|\bigl[|u|^2
		-|u(\xi_k)|^2\bigr]\circ R\biggr|\,dx\cr
		&\leqslant C_R^2\cdot [R(1)-R(0)]\,\varepsilon+4C_R\varepsilon\cdot
		\sum\limits_{k\in\Gamma}\int_{\zeta_k}^{\zeta_{k+1}}\kern -1mm
		f\,dx\cr
		&\leqslant\varepsilon\cdot\left[C_R^2\cdot[R(1)-R(0)]+4C_R\cdot
		\int_0^1 f\,dx\right].}
$$
Произвольность выбора величины $\varepsilon>0$ означает потому искомую возможность
представить меру $f$ в виде мультипликатора из пространства $\frak M_{R,U}$.
\endproof

\subsubsection\label{prop:2.1}
{\it Каков бы ни был мультипликатор $f\in\frak M_{R,U}$, оператор
$$
	S_{R,U}^*fS_{R,U}\colon W_{2,U}^1[R(0),R(1)]\to W_{2,U}^{-1}[R(0),R(1)]
$$
допускает представление в виде равномерного предела операторов умножения
на непрерывные функции.
}%
\proof
Ввиду характера определения пространства $\frak M_{R,U}$, достаточно прове\-рить
справедливость доказываемого предложения в случае $f\in L_1[0,1]$. В этом случае
ото\-бражение
$$
	u\mapsto\int_0^1 f\cdot\overline{(u\circ R)}\,dx
$$
является ограниченным функционалом на пространстве $W_2^1[R(0),R(1)]$, а потому
допускает запись в виде $u\mapsto\langle g,u\rangle$ с некоторой
$g\in W_2^{-1}[R(0),R(1)]$. Рассматривая последнюю обобщённую функцию в качестве
предела $g=\lim_{n\to\infty} g_n$ последовательности непрерывных функций, как раз
и получаем искомое представление оператора $S_{R,U}^*fS_{R,U}$.
\endproof

Отметим, что мультипликатор $S^*_{R,U}fS_{R,U}$, отвечающий суммируемой~--- и даже
непре\-рывной~--- функции $f$, может, вообще говоря, оказаться сингулярным.
Так, в слу\-чае, когда $R$~--- канторова лестница, а $f$~--- индикатор отрезка
$[2/5,3/5]$, мультипликатор $S^*_{R,U}fS_{R,U}$ представляет собой дельта-функцию
$(1/5)\cdot\bolddelta_{1/2}$.

\subsection\label{pt:1.2}
Зафиксируем пару самосопряжённых мультипликаторов $q,p\in\frak M_{R,U}$. Формальная
граничная задача
$$
	\displaylines{\rlap{$(\equation)$\label{eq:1}}\hfill
		-\left({y'\over r}\right)'+qy=\lambda py,\hfill\cr
		\rlap{$(\equation)$\label{eq:2}}\hfill
		(U-1)y^\vee+i(U+1)y^\wedge=0,\hfill}
$$
где через $y^\wedge$ и $y^\vee$, следуя [\cite{RH:2001}, \hbox{(7.50)}],
обозначены~--- также формальные~--- векторы граничных значений
$$
	y^\wedge\rightleftharpoons\pmatrix{y(0)\cr y(1)},\qquad
	y^\vee\rightleftharpoons\pmatrix{[y'/r](0)\cr -[y'/r](1)},
$$
допускает понимание в виде спектральной задачи для линейного операторного пучка
$T\colon\Bbb C\to{\cal B}(\frak D_{R,U},\frak D^*_{R,U})$, определённого правилом
$$
	\langle T(\lambda)y,z\rangle\equiv\int_{R(0)}^{R(1)} u'\overline{v'}\,dx+
		\langle [q-\lambda p]\cdot y,z\rangle+
		\langle Vu^\wedge,v^\wedge\rangle_{\Bbb C^2}.
	\leqno(\equation)
$$\label{eq:T}%
Здесь положено $u\rightleftharpoons S_{R,U}^{-1}y$ и $v\rightleftharpoons
S_{R,U}^{-1}z$, а символом $V\in\Bbb C^{2\times 2}$ обозначена самосопряжённая
матрица со свойством $(U-1)V=-i(U+1)$. В случае суммируемых $r>0$, $q$ и $p$
интегрированием по частям легко проверяется согласованность такого понимания
задачи \eqref{eq:1}, \eqref{eq:2} с обычным. Близкие к указанному~---
хотя формулируемые в иных терминах~--- понимания рассматриваемой задачи в частном
случае, когда коэффициенты $q$ и $p$ являются мерами или зарядами (с некоторыми
дополнительными ограничениями), могут быть найдены в работах [\cite{Fr:2011},
\cite{ET}]. Для случая равномерно положительной $r\in L_{\infty}[0,1]$ идентичная
рассмотренной трактовка задачи \eqref{eq:1}, \eqref{eq:2} введена в основанных
на идеях работ [\cite{NSh:1999}, \cite{SSh:1999}] работах [\cite{Vl:2004},
\cite{Vl:2009}].

Предложение \ref{prop:2.1} означает, что операторы $\hat T(\lambda)
\rightleftharpoons S_{R,U}^*T(\lambda)S_{R,U}$ допускают представление
$$
	\langle\hat T(\lambda)u,v\rangle\equiv\int_{R(0)}^{R(1)} u'\overline{v'}
		\,dx+\langle[\hat q-\lambda\hat p]\cdot u,v\rangle+
		\langle Vu^\wedge,v^\wedge\rangle_{\Bbb C^2},
$$
где $\hat q\rightleftharpoons S_{R,U}^*qS_{R,U}$ и $\hat p\rightleftharpoons
S_{R,U}^*pS_{R,U}$ суть вполне непрерывные мультипликаторы типа $W_{2,U}^1[R(0),
R(1)]\to W_{2,U}^{-1}[R(0),R(1)]$. Соответственно, вопрос об исследовании спектра
пучка $T$ вида \eqref{eq:T} может быть сведён к вопросу об исследовании спектра
граничной задачи
$$
	\displaylines{-u''+\hat qu=\lambda\hat pu,\cr
		(U-1)u^\vee+i(U+1)u^\wedge=0,}
$$
понимаемой в смысле работ [\cite{Vl:2004}, \cite{Vl:2009}]. Это наблюдение будет
широко использоваться нами в дальнейшем. В частности, из него вытекает
справедливость следующих предложений.

\subsubsection
{\it Если резольвентное множество пучка $T\colon\Bbb C\to{\cal B}(\frak D_{R,U},
\frak D^*_{R,U})$ вида~\eqref{eq:T} непусто, то спектр этого пучка является
дискретным.
}

\subsubsection\label{prop:PD}
{\it Пусть функции $\hat Q,\hat P\in L_2[R(0),R(1)]$ и числа $\hat Q(R(0))$,
$\hat Q(R(1))$, $\hat P(R(0))$ и $\hat P(R(1))$ таковы, что обобщённые функции
$\hat q,\hat p\in W_2^{-1}[R(0),R(1)]$ допускают представления
$$
	\leqalignno{\langle\hat q,u\rangle&\equiv-\int_{R(0)}^{R(1)}
		\hat Q\overline{u'}\,dx+\left.(\hat Q\overline{u})
		\right|_{R(0)}^{R(1)},&\eqlabel{eq:hatq}\cr
		\langle\hat p,u\rangle&\equiv-\int_{R(0)}^{R(1)}\hat P
		\overline{u'}\,dx+\left.(\hat P\overline{u})\right|_{R(0)}^{R(1)}.
	}
$$
Тогда соотношение $y\in\ker T(\lambda)$ равносильно существованию удовлетворяющего
равенству $y=u\circ R$ решения граничной задачи
$$
	\displaylines{\pmatrix{u\cr u^{[1]}}'=\pmatrix{\hat Q-\lambda\hat P&1\cr
		-[\hat Q-\lambda\hat P]^2&-\hat Q+\lambda\hat P}\cdot
		\pmatrix{u\cr u^{[1]}},\cr
		(U-1)\pmatrix{u^{[1]}(R(0))+[\hat Q(R(0))-\lambda\hat P(R(0))]\cdot
		u(R(0))\cr -u^{[1]}(R(1))-[\hat Q(R(1))-\lambda\hat P(R(1))]\cdot
		u(R(1))}+i(U+1)\pmatrix{u(R(0))\cr u(R(1))}=0.
	}
$$
В случае, когда при некотором $\lambda\in\Bbb C$ указанная граничная задача не имеет
нетривиальных решений, соответствующий оператор $T(\lambda)$ обладает ограниченным
обратным.
}

\subsection
Обозначим символом $A_p\rightleftharpoons\sup\{p,-p\}$ абсолютную величину
мультипликатора $p$ как самосопряжённого оператора класса ${\cal B}(\frak D_{R,U},
\frak D^*_{R,U})$. В частности, в случае неотрицательности мультипликатора $p$
справедливо равенство $A_p=p$. Введём в рассмотрение пространство $\frak H_{R,U,p}$,
получаемое пополнением пространства $\frak D_{R,U}$ по норме
$$
	\|y\|_{\frak H_{R,U,p}}\rightleftharpoons\sqrt{\langle A_py,y\rangle},
$$
и обозначим через $I\colon\frak D_{R,U}\to\frak H_{R,U,p}$ соответствующий
оператор плотного вложения. При этом оператор $I^*\colon\frak H_{R,U,p}\to
\frak D^*_{R,U}$ инъективен, а также может быть указана самосопряжённая инверсия
$J\colon\frak H_{R,U,p}\to\frak H_{R,U,p}$ со свойством $p=I^*JI$. Имеют место
следующие два факта.

\subsubsection\label{prop:3.1}
{\it Пусть при некотором $\xi\in\Bbb R$ оператор $T(\xi)$ вида \eqref{eq:T}
положителен. Тогда оператор $T^\bullet\rightleftharpoons (I^*J)^{-1}T(0)I^{-1}$
корректно определён и является \hbox{$J$-са}\-мо\-со\-пря\-жён\-ным.
}%
\proof
Для любых векторов $y,f\in\frak H_{R,U,p}$ разрешимость относительно
$z\in\frak D_{R,U}$ системы
$$
	Iz=y,\qquad T(0)z=I^*Jf
$$
равносильна разрешимости системы
$$
	Iz=y,\qquad T(\xi)z=I^*J(f-\xi y).
$$
Иначе говоря, в случае, когда ограниченный \hbox{$J$-са}\-мо\-со\-пря\-жён\-ный
оператор $I\cdot[T(\xi)]^{-1}I^*J$ имеет плотный образ, оператор $T^\bullet$
корректно определён и удовлетворяет равенству
$$
	T^\bullet=\xi+\bigl\{I\cdot[T(\xi)]^{-1}I^*J\bigr\}^{-1}.
$$
С целью установления плотности образа оператора $I\cdot[T(\xi)]^{-1}I^*J$ проверим
инъектив\-ность этого оператора. Пусть вектор $y\in\frak H_{R,U,p}$ удовлетворяет
равенству $I\cdot [T(\xi)]^{-1}I^*Jy\penalty 7500=0$. Тогда выполняются равенства
$$
	\eqalign{0&=\langle I\cdot[T(\xi)]^{-1}I^*Jy,Jy\rangle\cr
		&=\langle T(\xi)\cdot\{[T(\xi)]^{-1}I^*J\}y,
			[T(\xi)]^{-1}I^*Jy\rangle,
	}
$$
ввиду положительности оператора $T(\xi)$ означающие обращение в нуль векторов
$[T(\xi)]^{-1}I^*Jy$ и $I^*Jy$. Инъективность операторов $I^*$ и $J$ означает
теперь справедливость равенства $y=0$.
\endproof

\subsubsection
{\it При выполнении условий предложения \ref{prop:3.1} спектры оператора
$T^\bullet$ и пучка $T$ совпадают.
}

\proof
Непосредственно легко проверяется, что в случае принадлеж\-ности точки $\lambda\in
\Bbb C$ резольвентному множеству пучка $T$ оператор $I\cdot [T(\lambda)]^{-1}I^*J$
является обрат\-ным к оператору $T^\bullet-\lambda$. Соответственно, спектр
оператора $T^\bullet$ образован изолированными собственными значениями.
Для завершения доказательства теперь остаётся заметить, что при любом выборе числа
$\lambda\in\Bbb R$ и ненулевого вектора $y\in\ker T(\lambda)$ вектор
$Iy\in\frak H_{R,U,p}$ отличен от нуля и удовлетворяет равенству $T^\bullet Iy=
\lambda Iy$.
\endproof

Равенство $T^\bullet y=f$ может быть охарактеризовано в терминах
квазидифференциаль\-ных выражений. В частности, имеет место следующий факт.

\subsubsection\label{prop:AntiFr}
{\it Пусть мультипликатор $\hat q\in W_2^{-1}[R(0),R(1)]$ допускает представление
\eqref{eq:hatq}, а функция $f\in\frak H_{R,U,p}$ допускает представление
$$
	\langle f,JIS_{R,U}u\rangle\equiv-\int_{R(0)}^{R(1)}\hat F\overline{u'}\,dx
		+(\hat F\overline{u})\Bigr|_{R(0)}^{R(1)},
$$
где $\hat F\in L_2[R(0),R(1)]$ и $\hat F(R(0)),\hat F(R(1))\in\Bbb C$. Тогда
соотношение $T^\bullet y=f$ равносильно существованию удовлетворяющего равенству
$y=u\circ R$ решения граничной задачи
$$
	\displaylines{\pmatrix{u\cr u^{[1]}}'=\pmatrix{\hat Q&1\cr
		-\hat Q^2&-\hat Q}\cdot\pmatrix{u\cr u^{[1]}}+
		\pmatrix{-\hat F\cr \hat Q\hat F},\cr
		(U-1)\pmatrix{u^{[1]}(R(0))+\hat Q(R(0))u(R(0))-\hat F(R(0))\cr
		-u^{[1]}(R(1))-\hat Q(R(1))u(R(1))+\hat F(R(1))}
		+i(U+1)\pmatrix{u(R(0))\cr u(R(1))}=0.
	}
$$
}

Справедливость предложения \ref{prop:AntiFr} легко проверяется на основе
определения \eqref{eq:T} интегрирова\-нием по частям, аналогично стандартному
доказательству предложения \ref{prop:PD}. В случае, когда обобщённые функции $r$
и $p$ представляют собой непрерывные меры со свойством $\supp p\subseteq\supp r$,
а $q=0$, даваемое предложением \ref{prop:AntiFr} описание оператора $T^\bullet$
эквивалентно описанию операторов, рассматриваемых в работе~[\cite{Fr:2011}].


\section{Случай коэффициентов с самоподобными первообразными}
\subsection
В настоящем параграфе мы намерены изучить асимптотическое поведение спектра
задачи \eqref{eq:1}, \eqref{eq:2} в случае, когда $r$ есть вероятностная
непрерывная самоподобная мера, $q=0$, а $p$~--- обобщённая производная некоторой
самоподобной функции $P$. При этом мы будем опираться на указанную
в пункте~\ref{pt:1.2} возможность свед\'{е}ния рассматриваемой задачи к задаче
об обычной самоподобной струне, изучавшейся в работах [\cite{SV:1995}--%
\cite{VSh:2013}].

\subsection
Рассмотрим сначала случай, когда вероятностная мера $r$ есть обобщённая
производная неубывающей самоподобной функции $R\in C[0,1]$ с параметрами
$N>1$, $\{a_i\}_{i=1}^N$, $\{d_i\}_{i=1}^N$ и $\beta_i=\sum_{k<i}d_k$
(см. [\cite{VSh:2006:1}, \cite{Sh:2007}]). При этом естественным образом
предполагается выполнение равенств
$$
	\sum\limits_{i=1}^N a_i=\sum\limits_{i=1}^N d_i=1.
$$
Самоподобие обобщённой первообразной $P\in L_2([0,1];\,r)$ коэффициента $p$,
то есть функции со свойством
$$
	\langle p,u\circ R\rangle\equiv-\int_0^1 P\cdot\overline{u'\circ R}\,dR+
		\left.(P\overline{u})\right|_0^1,\leqno(\equation)
$$\label{eq:2.1}%
мы далее будем считать согласованным с самоподобием функции $R$ в том смысле,
что параметры $N$ и $\{a_i\}_{i=1}^N$ обеих функций являются общими. Для прочих
параметров самоподобия функции $P$ мы будем использовать обозначения
$\{d_i'\}_{i=1}^N$ и $\{\beta_i'\}_{i=1}^N$.

Стандартным образом (ср.~[\cite{VSh:2006:1}, \cite{Sh:2007}]) на основе принципа
сжимающих отображений устанав\-ливается справедливость следующего простого
предложения.

\subsubsection
{\it Функция $P\in L_2([0,1];\,r)$ с параметрами самоподобия $N$,
$\{a_i\}_{i=1}^N$, $\{d_i'\}_{i=1}^N$ и $\{\beta_i'\}_{i=1}^N$ существует
и однозначно определена, если выполнено неравенство
$$
	\sum\limits_{i=1}^N d_i\cdot|d_i'|^2<1.
$$
}

Построенная на основе самоподобной функции $P$ и произвольно фиксированной пары
"<граничных значений"> $P(0)$ и $P(1)$ по правилу \eqref{eq:2.1} обобщённая
функция $p$ очевидным образом принадлежит пространству $\frak M_{R,U}$
независимо от выбора матрицы граничных условий $U$. Соответствующая обобщённая
функция $\hat p\in W_2^{-1}[0,1]$ при этом, как легко выводится
из \eqref{eq:2.1}, допускает представление
$$
	\langle\hat p,u\rangle\equiv-\int_0^1\hat P\cdot\overline{u'}\,dx+
		\left.(P\overline{u})\right|_0^1,
$$
где $\hat P$~--- квадратично суммируемая самоподобная функция с параметрами
$M$, $\{d_{i_j}\}_{j=1}^M$, $\{d'_{i_j}\}_{j=1}^M$ и $\{\beta'_{i_j}\}_{j=1}^M$.
Через $M$ здесь обозначено число отличных от нуля элементов набора
$\{d_i\}_{i=1}^N$, а через $i_j$~--- нумерующие эти элементы индексы.

Утверждения об асимптотических свойствах спектра задачи \eqref{eq:1}, \eqref{eq:2}
с коэффициен\-тами рассматриваемого вида теперь немедленно выводятся из известных
[\cite{SV:1995}--\cite{VSh:2010}] спектральных свойств такого рода задач
с $r\equiv 1$.

\subsubsection [\cite{VSh:2006:1}, Теорема~4.2]
{\it Пусть набор $\{d_id'_i\}_{i=1}^N$ содержит не менее двух отличных от нуля
элементов, и пусть $\nu>0$~--- наибольший общий делитель конечных элементов набора
$\{-\ln (d_i\,|d'_i|)\}_{i=1}^N$, причём для некоторого индекса
$k\in\overline{1,N}$ выполняется одно из следующих условий:

\enumerate{9}{%
\item Справедливо неравенство $d'_k>0$, а величина $-\nu^{-1}\ln (d_k\,|d'_k|)$
нечётна.
\item Справедливо неравенство $d'_k<0$, а величина $-\nu^{-1}\ln (d_k\,|d'_k|)$
чётна.
}
\noindent
Пусть также величина $D>0$ удовлетворяет уравнению
$$
	\sum\limits_{i=1}^N (d_i\,|d'_i|)^D=1.\leqno(\equation)
$$\label{eq:D}%
Тогда считающие функции ${\bf N}_{\pm}$ собственных значений задачи \eqref{eq:1},
\eqref{eq:2} имеют при $\lambda\to\pm\infty$ асимптотики
$$
	{\bf N}_{\pm}(\lambda)= s_{\pm}(\ln |\lambda|)\cdot |\lambda|^D\cdot
		[1+o(1)],\leqno(\equation)
$$\label{eq:N}%
где $s_{\pm}$~--- непрерывные \hbox{$\nu$-пе}\-ри\-оди\-чес\-кие функции. В случае,
когда хотя бы при одном $k\in\overline{1,N}$ выполняется неравенство $d'_k<0$,
справедливо тождество $s_{-}(x)\equiv s_{+}(x)$.
}

\subsubsection [\cite{VSh:2006:2}, Теорема~4.1]
{\it Пусть набор $\{d_id'_i\}_{i=1}^N$ содержит не менее двух отличных от нуля
элементов, и пусть $\nu>0$~--- наибольший общий делитель конечных элементов набора
$\{-\ln (d_i\,|d'_i|)\}_{i=1}^N$. Пусть также для любого $k\in\overline{1,N}$
со свойством $d'_k>0$ величина $-\nu^{-1}\ln (d_k\,|d'_k|)$ чётна, а для любого
$k\in\overline{1,N}$ со свойством $d'_k<0$~--- нечётна. Тогда для считающих
функций ${\bf N}_{\pm}$ собственных значений задачи \eqref{eq:1}, \eqref{eq:2}
при $\lambda\to\pm\infty$ справедливы асимптотики \eqref{eq:N}, где $D>0$ есть
решение уравнения \eqref{eq:D}, а $s_{\pm}$~--- непрерывные \hbox{$2\nu$-пе}\-%
ри\-оди\-чес\-кие функции со свойством $s_{-}(x)\equiv s_{+}(x-\nu)$.
}

\subsubsection [\cite{VSh:2006:2}, Теорема~4.2]
{\it Пусть набор $\{d_id'_i\}_{i=1}^N$ содержит не менее двух отличных от нуля
элементов, причём конечные элементы набора $\{-\ln (d_i\,|d'_i|)\}_{i=1}^N$
не имеют общего делителя. Тогда для считающих функций ${\bf N}_{\pm}$ собственных
значений задачи \eqref{eq:1}, \eqref{eq:2} при $\lambda\to\pm\infty$ справедливы
асимптотики \eqref{eq:N}, где $D>0$ есть решение уравнения \eqref{eq:D},
а $s_{\pm}$~--- постоянные функции. В случае, когда хотя бы при одном $k\in
\overline{1,N}$ выполняется неравенство $d'_k<0$, справедливо равенство
$s_{-}=s_{+}$.
}

\subsubsection [\cite{VSh:2010}, Теоремы~4.1 и~4.2]
{\it Пусть набор $\{d_id'_i\}_{i=1}^N$ содержит ровно один отличный от нуля
элемент $d_md'_m>0$. Пусть также каждая из точек множества $\{R(\alpha_m),
R(\alpha_{m+1})\}\cap (0,1)$ является точкой разрыва функции $\hat P$.
Пусть, наконец, ${\rm Z}_{+}$ есть число точек $R(\alpha_i)\in (0,1)$, в которых
скачок функции $\hat P$ положителен, а ${\rm Z}_{-}$~--- число аналогичных точек,
в которых указанный скачок отрицателен. Тогда существуют такие вещественные числа
$\mu_{\pm,l}>0$, где $l\in\overline{1,{\rm Z}_{\pm}}$, что собственные значения
задачи \eqref{eq:1}, \eqref{eq:2} имеют при $k\to\infty$ асимптотики
$$
	\lambda_{\pm(l+k{\rm Z}_{\pm})}=\pm\mu_{\pm,l}\cdot (d_md'_m)^{-k}\cdot
		[1+o(1)].
$$
}

\subsubsection [\cite{VSh:2010}, Теорема~4.3]
{\it Пусть набор $\{d_id'_i\}_{i=1}^N$ содержит ровно один отличный от нуля
элемент $d_md'_m<0$. Пусть также каждая из точек множества $\{R(\alpha_m),
R(\alpha_{m+1})\}\cap (0,1)$ является точкой разрыва функции $\hat P$.
Пусть, наконец, ${\rm Z}_{+}$ есть число точек $R(\alpha_i)\in (0,1)$, в которых
скачок функции $\hat P$ положителен, а ${\rm Z}_{-}$~--- число аналогичных точек,
в которых указанный скачок отрицателен. Тогда существуют такие вещественные числа
$\mu_l>0$, где $l\in\overline{1,{\rm Z}_{+}+{\rm Z}_{-}}$, что собственные значения
задачи \eqref{eq:1}, \eqref{eq:2} имеют при $k\to\infty$ асимптотики
$$
	\eqalign{\lambda_{l+k({\rm Z}_{+}+{\rm Z}_{-})}&=\mu_l\cdot
		(-d_md'_m)^{-2k}\cdot [1+o(1)],\cr
		\lambda_{-(l+k{\rm Z}_{+}+(k+1){\rm Z}_{-})}&=-\mu_l\cdot
		(-d_md'_m)^{-2k-1}\cdot [1+o(1)].}
$$
}

\subsection
Пусть теперь обе функции $R$ и $P$ непрерывны, однако первая из них имеет несколько
более общий вид, чем ранее. А именно, пусть найдётся замкнутое подмно\-жество
$\Delta\subseteq [0,1]$, для которого независимо от выбора индекса
$i\in\overline{1,N}$ справедливо тождество
$$
	\int_{[\alpha_i,\alpha_i+a_ix]\cap\Delta}\,dR\equiv d_i\cdot R(x).
$$
Здесь введено обозначение $\alpha_i\rightleftharpoons\sum_{k<i}a_k$. Коэффициенты
$d_i$ при этом очевидным образом должны подчиняться условию
$$
	\sum_{i=1}^N d_i=\int_{\Delta}\,dR\leqslant 1.
$$
Коэффициенты самоподобия $N$ и $\{a_i\}_{i=1}^N$ функции $P$ по-прежнему
предполагаются совпада\-ющими с аналогичными коэффициентами "`почти самоподобия"'
функции $R$. Кроме того, мы будем предполагать, что для всякого $i\in\overline{1,N}$
со свойством $d_i'\neq 0$ выполняется также неравенство $d_i\neq 0$.

Ситуация из предыдущего пункта отвечает равенству $\Delta=[0,1]$.

При сделанных предположениях справедливо соотношение $\supp p\subseteq\supp r$,
означающее существование функции $\hat P\in C[0,1]$ со свойством $P=\hat P
\circ R$. Введём в рассмотрение функцию $\varphi:[0,1]\to [0,1]$ вида
$$
	\varphi(R(x))\equiv\int_{[0,x]\cap\Delta}\,dR.
$$
Тогда для всякого индекса $i\in\overline{1,N}$, удовлетворяющего неравенству
$d'_i\neq 0$, справедливо тождество
$$
	(\forall x\in [R(\alpha_i),R(\alpha_{i+1})])\qquad\hat P(x)=
		\beta'_i+d'_i\cdot\hat P\left({\varphi(x)-\varphi(R(\alpha_i))
		\over d_i}\right).
$$
Также имеет место следующий легко устанавливаемый на основе метода расщепления
(по точкам излома функции $\varphi$) факт.

\subsubsection\label{prop:2.3.1}
{\it Пусть $\Delta$ есть объединение конечного числа попарно непересекающихся
отрезков. Тогда при всяком $i\in\overline{1,N}$ со свойством $d'_i\neq 0$ индекс
инерции заданной на пространстве $\Wo_2^1[R(\alpha_i),R(\alpha_{i+1})]$
квадратичной формы
$$
	y\mapsto\int_{R(\alpha_i)}^{R(\alpha_{i+1})}\left[|y'(x)|^2+
		\lambda\hat P\left({\varphi(x)-\varphi(R(\alpha_i))
		\over d_i}\right)\cdot (|y|^2)'(x)\right]\,dx
$$
имеет при $\lambda\to+\infty$ асимптотику ${\bf N}(d_i\lambda)+O(1)$, где через
${\bf N}(\lambda)$ обозначен индекс инерции заданной на пространстве $\Wo_2^1[0,1]$
квадратичной формы
$$
	y\mapsto\int_0^1\bigl[|y'|^2+\lambda\hat P\cdot (|y|^2)'\bigr]\,dx.
$$
}

\noindent
Повторяя с учётом предложения \ref{prop:2.3.1} рассуждения из доказательств
предложений [\cite{VSh:2006:1}, Теорема~4.2] и [\cite{VSh:2006:2}, Теоремы
4.1 и~4.2], приходим к следующему факту.

\subsubsection
{\it Пусть $\Delta$ есть объединение конечного числа попарно непересекающихся
отрезков. Тогда для считающих функций ${\bf N}_{\pm}$ собственных значений
задачи \eqref{eq:1}, \eqref{eq:2} при $\lambda\to\pm\infty$ справедливы асимптотики
\eqref{eq:N}, где $s_{\pm}$~--- непрерывные периодические функции, а $D$~---
положительное решение уравнения \eqref{eq:D}. В случае, когда конечные элементы
набора $\{-\ln (d_i\,|d'_i|)\}_{i=1}^N$ не имеют общего делителя, функции $s_{\pm}$
являются постоянными.
}

\subsection
Результаты, касающиеся описанного в предыдущем пункте случая с имеющим общий вид
замкнутым множеством $\Delta\subseteq [0,1]$, заявлены в работе
[\cite{Fr:2011}]. Предложенные в указанной работе доказательства, однако,
являются ошибочными, поскольку опираются на следующее легко опровергаемое
утверждение.

\subsubsection\label{lemma} [\cite{Fr:2011}, Lemma~3.1]
{\it Считающая функция $\bf N$ собственных значений задачи \eqref{eq:1},
\eqref{eq:2} с "`почти самоподобной"' $R$, неубывающей непрерывной самоподобной $P$
и матрицей граничных условий Неймана $U=-1$ при всех $\lambda\geqslant 0$
удовлетворяет неравенству
$$
	{\bf N}(\lambda)\leqslant\sum\limits_{i=1}^N{\bf N}(d_id'_i\cdot\lambda).
$$
}

\topinsert
\epsfxsize=0.9\hsize
\hskip 0pt\hbox to \epsfxsize{\epsfbox[-8 -9 420 126]{Iter.eps}}\par
\hbox to \hsize{\hfill Рис.~1.\hfill}
\endinsert
Действительно, рассмотрим в качестве $P$ канторову лестницу с параметрами
самоподобия $N=3$, $a_1=a_2=a_3=1/3$, $d'_1=d'_3=1/2$, $d'_2=0$, $\beta'_1=0$,
$\beta'_2=\beta'_3=1/2$, а в качестве $R$~--- результат $P_3$ третьего шага
итеративного построения той же канторовой лестницы, отправляющегося от функции
$P_0(x)\equiv x$ (см. Рис.~1). Параметры "`почти самоподобия"' функции $P_3$
имеют вид $d_1=d_3=1/3$, $d_2=0$. Соответственно (см.~Таб.~1), соотношения
$$
	{\bf N}(510)=9>3+1+3\geqslant{\bf N}(85)+{\bf N}(0)+{\bf N}(85)
$$
опровергают утверждение \ref{lemma}.
\topinsert
\rightskip=0pt plus 1fill\leftskip=0pt plus 1fill\leavevmode
\vbox{\offinterlineskip\hrule
	\halign{&\vrule#&\strut\quad\hfill#\hfill\quad\cr
	height 1pt&\omit&&\omit&&\omit&&\omit&\cr
	&$n$\hfil&&$\lambda_n$\hfil&&$n$\hfil&&$\lambda_n$\hfil&\cr
	height 1pt&\omit&&\omit&&\omit&&\omit&\cr
	\noalign{\hrule}
	height 1pt&\omit&&\omit&&\omit&&\omit&\cr
	&$1$&&$0{,}000\,({\pm 10^{-3}})$&&$6$&&$235{,}7\,({\pm 10^{-1}})$&\cr
	&$2$&&$9{,}857\,({\pm 10^{-3}})$&&$7$&&$327{,}3\,({\pm 10^{-1}})$&\cr
	&$3$&&$39{,}27\,({\pm 10^{-2}})$&&$8$&&$413{,}4\,({\pm 10^{-1}})$&\cr
	&$4$&&$87{,}69\,({\pm 10^{-2}})$&&$9$&&$454{,}2\,({\pm 10^{-1}})$&\cr
	&$5$&&$154{,}0\,({\pm 10^{-1}})$&&$10$&&$982{,}3\,({\pm 10^{-1}})$&\cr
	height 1pt&\omit&&\omit&&\omit&&\omit&\cr}
	\hrule
}
\par\penalty 10000\medskip Таб.~1.
\endinsert

Таким образом, проблему описания спектральных асимптотик задачи \eqref{eq:1},
\eqref{eq:2} для произвольных согласованных "`почти самоподобных"' пар мер следует,
по всей вероятности, считать открытой.

\penalty 10000\vskip 1ex
Автор выражает признательность проф.~А.$\,$И.~Назарову и проф.~И.$\,$А.~Шейпаку
за ценные обсуждения.


\vskip 0.4cm
\eightpoint\rm
{\leftskip 0cm\rightskip 0cm plus 1fill\parindent 0cm
\bf Литература\par\penalty 20000}\vskip 0.4cm\penalty 20000
\bibitem{RH:2001} Ф.$\,$С.~Рофе--Бекетов, А.$\,$М.~Холькин. {\it Спектральный
анализ дифференциальных операторов. Связь спек\-тральных и осцилляционных свойств}.
Мариуполь, 2001.

\bibitem{Fr:2011} U.~Freiberg. {\it Refinement of the spectral asymptotics
of generalized Krein Feller operators}~// Forum Math.~--- 2011.~--- V.~23.~---
P.~427--445.

\bibitem{ET} J.~Eckhardt, G.~Teschl. {\it Sturm--Liouville operators with
measure-valued coefficients}~// Journal d'analyse math\'ematique.~--- 2013.~---
V.~120, \No~1.~--- P.~151--224.

\bibitem{NSh:1999} М.$\,$И.~Нейман--заде, А.$\,$А.~Шкаликов. {\it Операторы
Шрёдингера с сингулярными потенциалами из прост\-ранств мультипликаторов}~//
Матем. заметки.~--- 1999.~--- Т.~66, \No~5.~--- С.~723--733.

\bibitem{SSh:1999} А.$\,$М.~Савчук, А.$\,$А.~Шкаликов. {\it Операторы
Штурма--Лиувилля с сингулярными потенциа\-лами}~// Матем. заметки.~--- 1999.~---
Т.~66, \No~6.~--- С.~897--912.

\bibitem{Vl:2004} А.$\,$А.~Владимиров. {\it О сходимости последовательностей
обыкновенных дифференциальных операторов}~// Матем. заметки.~--- 2004.~---
Т.~75, \No~6.~--- С.~941--943.

\bibitem{Vl:2009} А.$\,$А.~Владимиров. {\it К осцилляционной теории задачи
Штурма--Лиувилля с сингулярными коэффициен\-тами}~// Ж.~выч. матем.
и матем.~физ.~--- 2009.~--- Т.~49, \No~9.~--- С.~1609--1621.

\bibitem{SV:1995} M.~Solomyak, E.~Verbitsky. {\it On a spectral problem
related to self-similar measures}~// Bull. London Math. Soc.~--- 1995.~---
V.~27.~--- P.~242--248.

\bibitem{VSh:2006:1} А.$\,$А.~Владимиров, И.$\,$А.~Шейпак. {\it Самоподобные
функции в пространстве $L_2[0,1]$ и задача Штурма--Лиувилля с сингулярным
индефинитным весом}~// Матем. сб.~--- 2006.~--- Т.~197, \No~11.~--- С.~13--30.

\bibitem{VSh:2006:2} А.$\,$А.~Владимиров, И.$\,$А.~Шейпак. {\it Индефинитная
задача Штурма--Лиувилля для некоторых классов самоподобных сингулярных весов}~//
Труды МИАН.~--- 2006.~--- Т.~255.~--- С.~88--98.

\bibitem{VSh:2010} А.$\,$А.~Владимиров, И.$\,$А.~Шейпак. {\it Асимптотика
собственных значений задачи Штурма--Лиувилля с дискретным самоподобным
весом}~// Матем. заметки.~--- 2010.~--- Т.~88, \No~5.~--- С.~662--672.

\bibitem{VSh:2013} А.$\,$А.~Владимиров, И.$\,$А.~Шейпак. {\it О задаче Неймана
для уравнения Штурма--Лиувилля с самоподобным весом канторовского типа}~//
Функц. анализ и его прилож.~--- 2013.~--- Т.~47, \No~4.~--- С.~18--29.

\bibitem{Sh:2007} И.$\,$А.~Шейпак. {\it О конструкции и некоторых свойствах
самоподобных функций в простран\-ствах $L_p[0,1]$}~// Матем. заметки.~--- 2007.~---
Т.~81, \No~6.~--- С.~924--938.
\bye